\documentclass[reqno,11pt]{amsart}
\usepackage{amsthm}
\usepackage{amsmath,mathrsfs}
\usepackage{color}
\usepackage{amssymb}
\usepackage{hyperref}
\usepackage{enumerate}
\usepackage{latexsym,array}
\usepackage{amsfonts}
\usepackage{shadow}
\newcommand{\ran}{{\rm ran}\,}

\newtheorem{Pa}{Paper}[section]
\newtheorem{Tm}[Pa]{{\bf Theorem}}
\newtheorem{La}[Pa]{{\bf Lemma}}
\newtheorem{Dn}[Pa]{{\bf Definition}}
\newtheorem{Cy}[Pa]{{\bf Corollary}}

\newtheorem{Rk}[Pa]{{\bf Remark}}
\newtheorem{Pn}[Pa]{{\bf Proposition}}

\newcommand{\vv}{\varphi}
\newcommand{\e}{\varepsilon}

\def\hh{\mathbb{H}}
\author[D. Alpay]{Daniel Alpay}
\address{(DA) Department of Mathematics\\
Ben-Gurion University of the Negev\\
Beer-Sheva 84105 Israel} \email{dany@math.bgu.ac.il}

\author[F. Colombo]{Fabrizio Colombo}
\address{(FC) Politecnico di
Milano\\Dipartimento di Matematica\\Via E. Bonardi, 9\\20133
Milano, Italy}
\email{fabrizio.colombo@polimi.it}

\author[D. P. Kimsey]{David P. Kimsey}
\address{(D.P.K.) Department of Mathematics\\
Ben-Gurion University of the Negev\\
Beer-Sheva 84105 Israel}
\email{dpkimsey@gmail.com}
\author[I. Sabadini]{Irene Sabadini}
\address{(IS) Politecnico di
Milano\\Dipartimento di Matematica\\Via E. Bonardi, 9\\20133
Milano, Italy}
\email{irene.sabadini@polimi.it}
\title[An extension of Herglotz's theorem to the quaternions]
{An extension of Herglotz's theorem to the quaternions} \oddsidemargin
0.2in \evensidemargin 0.2in \topmargin -0.5in \textwidth
15.5truecm \textheight 23truecm

\def\R{\mathbb R}

\def\va{\varphi}

\def\(s){\mathscr S(\R\times\R)}



\thanks{D. Alpay thanks the Earl Katz family for endowing the chair
which supported his research. D.P. Kimsey gratefully acknowledges the support of a Kreitman postdoctoral fellowship. F. Colombo and I. Sabadini
acknowledge the Center for Advanced Studies of the Mathematical
Department of the Ben-Gurion University of the Negev for the
support and the kind hospitality during the period in which part
of this paper has been written.}
\begin{document}
\maketitle
\tableofcontents
\parindent 0cm
\begin{abstract}
A classical theorem of Herglotz states that a
function $n\mapsto r(n)$ from $\mathbb Z$ into $\mathbb
C^{s\times s}$ is positive definite if and only there exists a
$\mathbb C^{s\times s}$-valued positive measure $\mu$ on $[0,2\pi]$ such
that $r(n)=\int_0^{2\pi}e^{int}d\mu(t)$for $n\in \mathbb Z$.
We prove a quaternionic analogue of this result when the function is allowed to have a number of negative squares.
A key tool in the argument is the theory of slice hyperholomorphic functions, and the representation of such functions which
have a positive real part in the unit ball of the quaternions. We study in great detail the case of positive definite functions.
 \end{abstract}

\parindent 0cm

\section{Introduction}
\setcounter{equation}{0}
The main purpose of this paper is to prove a version of a theorem of Herglotz on
positive functions in the quaternionic and indefinite setting. To
set the framework we first recall some definitions and results
pertaining to the complex numbers setting. A function $n\mapsto
r(n)$ from $\mathbb Z$ into $\mathbb C^{s\times s}$ is called
positive definite if the associated function (also called kernel)
$K(n-m)$ is positive definite on $\mathbb Z$. This means that for
every choice of $N\in\mathbb N$ and $n_1,\ldots, n_N\in\mathbb
Z$, the $N\times N$ block matrix with $(j,\ell)$ block entry
equal to the matrix $r(n_j-n_\ell)$ is non-negative, that is, all
the block Toeplitz matrices
\begin{equation}
\label{eq:toeplitz}
\mathbb T_N\stackrel{\rm def.}{=}\begin{pmatrix}r(0)&r(1)&\cdots &r(N)\\
                                         r(-1)&r(0)&\cdots &r(N-1)\\
                                               & & &\\
                                               & & &\\
                                         r(-N)&r(1-N)&\cdots& r(0)\end{pmatrix}
\end{equation}
are non-negative. We will use the notation $\mathbb T_N \succeq 0$.
The positivity implies in particular that
$r(-n)=r(n)^*$, where $r(n)^*$ denotes the adjoint of $r(n)$.\smallskip

A result of Herglotz, also known as Bochner's theorem, asserts
that:

\begin{Tm}	
\label{Tm:Oct27yt1}
The function $n \mapsto r(n)$ from $\mathbb Z$ into $\mathbb C^{s \times s}$ is positive definite if and only if there exists a unique positive $\mathbb C^{s \times s}$-valued measure $\mu$ on $[0, 2\pi ]$ such that
\begin{equation}
\label{eq:Oct27nv1}
r(n) = \int_0^{2\pi} e^{i n t} d\mu(t), \quad n \in \mathbb Z.
\end{equation}
\end{Tm}
See for instance \cite[p. 38]{MR0422992}, \cite{herglotz} and the
discussion in \cite[p. 19]{MR1038803}. We note that
\eqref{eq:Oct27nv1} can be rewritten as
\begin{equation}
r(n)=C^*U^nC,\quad n\in \mathbb Z,
\end{equation}
where $U$ denotes the unitary operator of multiplication by
$e^{it}$ in $\mathbf L_2([0,2\pi],d\mu)$, and $C$ denotes the
operator defined by $C\xi=\xi$ from $\mathbb C^s$ into $\mathbf
L_2([0,2\pi],d\mu)$.\\

A result of Carath{\'e}odory \cite{Caratheodory}, which is related to Theorem \ref{Tm:Oct27yt1} asserts that:

\begin{Tm}	
\label{Tm:Jan10u1}
If the function $n \mapsto r(n)$ from ${ -N, \ldots, N}$ into $\mathbb C^{s \times s}$ is positive definite, i.e., $\mathbb{T}_N \succeq 0$, then exists a function $n \mapsto \tilde{r}(n)$ from $\mathbb{Z}$ to $\mathbb{C}^{s \times s}$ which is positive definite and satisfies
\begin{equation}
\label{eq:Jan10y1}
r(n) = \tilde{r}(n), \quad n \in \{-N, \ldots, N\}.
\end{equation}
\end{Tm}
\smallskip

In commutative harmonic analysis, Theorem \ref{Tm:Oct27yt1} is a
special case of a general result of Weil \cite{MR0005741} on the
representation of positive definite functions on a group in terms
of the characters of the group. See for instance \cite[(22.7.10),
p
65]{dieudonne-6} or \cite[Theorem 5.4.3, p. 65]{MR2640609}.\\

A key result in one of the proofs (see for instance \cite[pp.
148-149]{MR0422992}) of Theorem \ref{Tm:Oct27yt1} is Herglotz's
representation theorem, which states that a $\mathbb C^{s\times
s}$-valued function $\va$ is analytic and with a real positive
part in the open unit disk $\mathbb D$ if and only if it can be
written as
\begin{equation}
\label{vaher}
\va(z)=\int_0^{2\pi}\frac{e^{it}+z}{e^{it}-z}d\mu(t)+ia,
\end{equation}
where $d\mu$ is as in Theorem \ref{Tm:Oct27yt1} and $a\in\mathbb
C^{s\times s}$ satisfies $a+a^*=0$. There are a number of ways to
prove \eqref{vaher}. It can be obtained from Cauchy's formula and
from the weak-$*$ compactness of the family of finite variation measures on
$[0,2\pi]$; see for instance the discussion in \cite[p.
207]{CAPB}.\\

Krein extended the notion of positive definite functions to the
notion of functions having a number of negative squares; see
\cite{ikl}. We first recall the definition of this notion in the
present setting:

\begin{Dn}
\label{nsq} The function $n\mapsto r(n)$ from $\mathbb Z$ into
$\mathbb C^{s\times s}$ satisfying $r(n)=r(-n)^*$ has a finite
number of negative squares, say $\kappa$, if by definition the
function $K(n,m)=r(n-m)$ has $\kappa$ negative squares, that is,
if all the block Toeplitz matrices $\mathbb T_N$  defined in
\eqref{eq:toeplitz} (which are Hermitian since $r(n)=r(-n)^*$)
have at most $\kappa$ strictly negative eigenvalues and exactly
$\kappa$ strictly negative eigenvalues for some choice of $N$ and
$n_1,\ldots n_N$.
\end{Dn}
Theorem \ref{Tm:Oct27yt1} was extended, in the scalar case, by
Iohvidov \cite{MR728210}  to case where the function $K(n,m)$ has
a finite number of negative squares. Formula \eqref{eq:Oct27nv1} is
then replaced by a more involved expression. More precisely, he
obtained the following extension of \eqref{eq:Oct27nv1} (there is a
minus sign with respect to \cite{MR728210} and \cite{MR0107821}
because they work there with positive squares rather than
negative squares):
\begin{align}
r(n)  =& \; \int_{0}^{2\pi}\frac{e^{int}-S_n(t)}{\prod_{k=1}^u \left(
\sin\left(\frac{t-\vv_k}{2} \right)^{2\rho_k} \right)}d\mu(t) \nonumber \\
& \; \; -\left(\sum_{j=1}^rQ_j(in)\lambda_j^n+\overline{Q_j(in)\lambda_j}^{-n}
+\sum_{k=1}^uR_k(in)e^{in\vv_k}\right). \label{klformula}
\end{align}
In this expression, the $\lambda_j$ are of modulus strictly
bigger than $1$, the $Q_j$ and $R_j$ are polynomials and $S_n$ is
a regularizing correction. These terms follow from the structure
of a contraction in a Pontryagin space, and in particular from
the fact that such an operator has always a strictly negative
invariant subspace, on which it is one-to-one. See \cite[(20.2), p. 319]{MR0107821}, where Iohvidov and Krein prove that such
a representation is unique.\\

In this paper we shall prove in particular a quaternionic analogue of Theorem \ref{Tm:Oct27yt1}, where
$\mathbb{H}$ denote the quaternions:

\begin{Tm}
Let $(r(n))_{n\in\mathbb Z}$ be a sequence of $s\times s$ matrices
with quaternionic entries. Then:\\
$(1)$ The function $K(n,m)=r(n-m)$ has a finite number of
negative squares $\kappa$ if and only if there exists a right
quaternionic Pontryagin space $\mathcal P$, a unitary operator
$U\in\mathbf L(\mathcal P)$ and a linear operator $C\in\mathbf
L(\mathbb H^s,\mathcal P)$ such that
\begin{equation}
\label{voltaire}
r(n)=C^*U^nC,\quad n\in\mathbb Z.
\end{equation}
$(2)$ Assume that
\begin{equation}
\label{chatelet}
\bigcup_{n\in\mathbb Z}{\rm ran}\, U^nC
\end{equation}
is dense in $\mathcal P$ where {\rm ran} $U^nC$ denotes the range of $U^nC$. Then, the realization \eqref{voltaire} is unique up to a
unitary map.
\label{montmartre}
\end{Tm}
Some remarks:\\
$(1)$ The sufficiency of condition \eqref{voltaire} follows from
the inner product representation
\[
c^*r(n-m)d=c^*C^*U^{n-m}Cd=\langle\,
U^{-m}Cd\,,\,U^{-n}Cc\rangle_{\mathcal P},\quad n,m\in\mathbb
Z,\quad c,d\in\mathbb H^s.
\]
The proof of the necessity is done using the theory of slice
hyperholomorphic functions. We use in particular a representation theorem from
\cite{acs2} for functions $\varphi$ slice-hyperholomorphic in some open
subset of the unit ball and with a certain associated kernel $K_\varphi(p,q)$ (defined by \eqref{opera} below) having a finite number of negative squares  there. We also note that the arguments in \cite{acs2} rely on the theory of linear relations
in Pontryagin spaces.\smallskip

$(2)$ The more precise integral representation of Iohvidov and
Krein relies on the theory of unitary operators in Pontryagin
spaces. Such results are still lacking in the setting of
quaternionic Pontryagin spaces.\smallskip

$(3)$ We also consider the positive definite case. There, the lack of a properly
established spectral theorem for unitary operators in quaternionic
Hilbert spaces prevents to get a direct counterpart of the integral representation
\eqref{eq:Oct27nv1}.\\

The outline of the paper is as follows. The paper consists of seven
sections, besides the introduction. In Section \ref{sec:two} we
review some results from the theory of slice hyperholomorphic
functions.
 Some definitions and results on quaternionic
Pontryagin spaces are recalled in Section \ref{sec:3} as well as Herglotz-type theorem for matrix valued functions. The proof of
the necessity and uniqueness in Theorem \ref{montmartre} is done in Section
\ref{sec:three}. Section \ref{Bochner} deals with the analogue of Herglotz's theorem in the quaternionic setting.
In Section \ref{carath:sec} we prove a quaternionic analogue of Theorem \ref{Tm:Jan10u1}.
Section
  7 contains the characterization  of quaternionic, bounded, Hermitian sequences of matrices with $\kappa$ negative squares. It uses results proved in Section  \ref{Bochner}.
In Section \ref{sec:six} we prove an
Herglotz  representation theorem for scalar valued functions slice hyperholomophic
in the unit ball of the quaternions, and with a positive real part there.

\section{Slice hyperholomorphic functions}
\setcounter{equation}{0}
\label{sec:two}

The  kernels we will use in this paper are slice hyperholomorphic, so we recall their definition.
For more details and the proofs of the results in this section see \cite{MR2752913}.
\\
The imaginary units  in $\mathbb{H}$ are denoted by $i$, $j$ and $k$,
 and an element in $\mathbb{H}$ is of the form $p=x_0+ix_1+jx_2+kx_3$, for $x_\ell\in \mathbb{R}$.
The real part, the imaginary part and conjugate of $p$ are defined as
${\rm Re}(p)=x_0$,  ${\rm Im}(p)=i x_1 +j x_2 +k x_3$ and
by $\bar p=x_0-i x_1-j x_2-k x_3$, respectively.
\\
The unit sphere of purely imaginary  quaternions $\mathbb{S}$ is defined by
$$
\mathbb{S}=\{q=ix_1+jx_2+kx_3\ {\rm such \ that}\
x_1^2+x_2^2+x_3^2=1\}.
$$
Note that if $I\in\mathbb{S}$, then
$I^2=-1$; for this reason the elements of $\mathbb{S}$ are also called
imaginary units. Note that $\mathbb{S}$ is a 2-dimensional sphere in $\mathbb{R}^4$.
Given a nonreal quaternion $p=x_0+{\rm Im} (p)=x_0+I |{\rm Im} (p)|$, $I={\rm Im} (p)/|{\rm Im} (p)|\in\mathbb{S}$, we can associate to it the 2-dimensional sphere defined by
$$
[p]=\{x_0+I|{\rm Im} (p)|\ :\ I\in\mathbb{S}\}.
$$
We will denote an element in the complex plane $\mathbb{C}_I:=\mathbb{R}+I\mathbb{R}$ by $x+Iy$.
\begin{Dn}[Slice hyperholomorphic functions]\label{regularity} Let $\Omega$ be an open set in
$\mathbb{H}$ and let $f:\Omega \to \mathbb{H}$ be a real differentiable function.
Denote by  $f_I$ the restriction of $f$  to the complex plane $\mathbb{C}_I$.
\\
We say that  $f$  is (left) slice hyperholomorphic  (or (left) slice regular) if, on $\Omega \cap \mathbb{C}_I$, $f_I$ satisfies
$$\frac{1}{2}\left(\frac{\partial}{\partial x}
+I\frac{\partial}{\partial y}\right)f_I(x+Iy)=0,$$
for all $I \in
\mathbb{S}$.
\\
We say that  $f$ is
right slice hyperholomorphic (or right slice regular)  if, on $\Omega \cap \mathbb{C}_I$, $f_I$ satisfies
$$\frac{1}{2}\left(\frac{\partial}{\partial x}f_I(x+Iy)
+\frac{\partial}{\partial y}f_I(x+Iy)I\right)=0,$$
for all $I \in
\mathbb{S}$.
\end{Dn}
An immediate consequence of the definition of slice regularity is that
the monomial $p^na$, with $a\in\mathbb{H}$, is
left slice regular, so power series with quaternionic coefficients written on the right are
left slice regular where they converge.
As one can easily verify only power series with center at real points are slice regular.

We introduce a class of domains, which includes the balls with center at a real point, on which slice regular functions have good properties.
\begin{Dn}[Axially symmetric domain]
Let $U \subseteq \mathbb{H}$. We say that $U$ is
\textnormal{axially symmetric} if, for all $x+Iy \in U$, the whole
2-sphere $[x+Iy]$ is contained in $U$.
\end{Dn}

\begin{Dn}[Slice domain]
Let $U \subseteq \mathbb{H}$ be a domain in $\mathbb{H}$. We say that $U$ is a
\textnormal{slice domain (s-domain for short)} if $U \cap \mathbb{R}$ is non empty and
if $ U\cap\mathbb{C}_I$ is a domain in $\mathbb{C}_I$ for all $I \in \mathbb{S}$.
\end{Dn}

\begin{La}[Splitting Lemma]\label{splitting}
Let $\Omega$ be an  s-domain in $\mathbb{H}$.
If $f:\Omega \to \mathbb{H}$ is left slice hyperholomorphic, then for
every $I \in \mathbb{S}$, and every $J\in\mathbb{S}$,
perpendicular to $I$, there are two holomorphic functions
$F,G:\Omega \cap \mathbb{C}_I \to \mathbb{C}_I$ such that for any $z=x+Iy$, it is
$$f_I(z)=F(z)+G(z)J.$$
\end{La}
Note that the decomposition given in the Splitting Lemma is highly non-canonical. In fact, for any $I\in\mathbb S$ there is an infinite number of choices of $J\in\mathbb S$ orthogonal to it.
\begin{Tm}[Representation Formula] Let $\Omega$ be an axially symmetric s-domain $\Omega \subseteq  \mathbb{H}$ and
let $f:\Omega \to \mathbb{H}$ be a slice hyperholomorphic function on $\Omega$.  Then the following equality holds for all $x+yI, x\pm Jy \in \Omega $:
\begin{equation}\label{distribution}
f(x+Iy) =\frac{1}{2}\Big[   f(x+Jy)+f(x-Jy)\Big] +I\frac{1}{2}\Big[ J[f(x-Jy)-f(x+Jy)]\Big].
\end{equation}
\end{Tm}

\section{Quaternionic Pontryagin spaces}
\setcounter{equation}{0}
\label{sec:3}

A Hermitian form on a right quaternionic vector space $\mathcal P$
is an $\mathbb H$-valued map $[\cdot, \cdot]$ defined on $\mathcal
P\times \mathcal P$ and such that
\[
\begin{split}
[a,b]&=\overline{[b,a]}\\
[ap,bq]&=\overline{q}[a,b]p,\quad\forall a,b\in\mathcal P\,\,{\rm
and}\,\,p,q\in\mathbb H.
\end{split}
\]
$\mathcal P$ is called a (right quaternionic) Pontryagin space if
it can written as
\begin{equation}
\label{fundamental}
\mathcal P=\mathcal
P_+\stackrel{\cdot}{[+]}\mathcal P_-,
\end{equation}
where:\\
$(a)$ The space $\mathcal P_+$ endowed with the form  $[\cdot,
\cdot]$ is a Hilbert space.\\
$(b)$ The space $\mathcal P_-$ endowed with the form  $-[\cdot,
\cdot]$ is a finite dimensional Hilbert space.\\
$(c)$ The sum is direct and orthogonal, meaning that $\mathcal
P_+\cap\mathcal P_-=\left\{0\right\}$ and
\[
[a,b]=0,\quad\forall (a,b)\in\mathcal P_+\times \mathcal P_-.
\]
The decomposition \eqref{fundamental} is called a fundamental
decomposition. It not unique unless one of the components reduces
to $\left\{0\right\}$. The dimension of $\mathcal P_-$ is the same
for all fundamental decompositions, and is called the index of the
Pontryagin space. The space $\mathcal P$ endowed with the form
\begin{equation}
\label{hilbertinner}
\langle a,b\rangle=[a_+,b_+]-[a_-,b_-]
\end{equation}
with $a_\pm$ and $b_\pm\in\mathcal P_\pm$ is a Hilbert space. The
inner product \eqref{hilbertinner} depends on the given
fundamental decomposition, but all the associated norms are
equivalent,and hence define the same topology. We refer to
\cite{as3} for a proof of these facts.  We refer to
\cite{as3,2013arXiv1303.1076A} for more details on quaternionic
Pontryagin spaces and to\cite{azih,bognar,ikl} for the theory of Pontryagin
spaces in the complex case. A reproducing kernel Pontryagin space will be a
Pontryagin space of functions for which the point evaluations are bounded. The definition of negative squares makes sense in the
quaternionic setting since an Hermitian quaternionic matrix $H$
is diagonalizable: It can be written as $T=UDU^*$, where $U$ is
unitary and $D$ is unique and with real entries. The number of
strictly negative eigenvalues of $T$ is exactly the number of
strictly negative elements of $D$. See \cite{MR97h:15020}. The
one-to-one correspondence between reproducing kernel Pontryagin
spaces and functions with a finite number of negative squares,
proved in the classical case by \cite{schwartz,Sorjonen73},
extends to the Pontryagin space setting, see \cite{as3}.

\begin{Dn}
An $\mathbb{H}^{s\times s}$-valued function $\varphi$ slice hyperholomorphic in a neighborhood $\mathcal{V}$ of the
origin is called a generalized Carath\'eodory function if the kernel
$$
k_{\varphi}(p,q)=\Sigma_{\ell=0}^\infty p^\ell(\varphi(p)+\overline{\varphi(q)})\overline{q}^\ell
$$
has a finite number of negative squares in $\mathcal{V}$.
\end{Dn}
The following result is Theorem 10.2 in \cite{acs2}.
\begin{Tm}
A $\mathbb{H}^{s\times s}$-valued function $\varphi$
is a generalized Carath\'eodory function if and only if it can be written as
\begin{equation}\label{REPPP}
\varphi(p)=\frac{1}{2} C\star (I_{\mathcal P}+pV)\star (I_{\mathcal P}-pV)^{-\star} C^*J+\frac{\varphi(0)-\varphi(0)^*}{2}
\end{equation}
where ${\mathcal P}$ is a right quaternionic Pontryagin space of index $\kappa$, $V$ is a coisometry in ${\mathcal P}$,
and $C$ is a bounded operator from ${\mathcal P}$ to $\mathbb{H}^N$, and the pair $(C, A)$ is observable.
\end{Tm}
\begin{Rk}{\rm
When $\kappa=0$ the representation (\ref{REPPP}) is the counterpart of Herglotz representation theorem for functions slice hyperholomorphic in the open unit ball and with a positive real part. In the last section we shall discuss a scalar version of this result.
}
\end{Rk}

\smallskip

\section{Proof of the necessity and uniqueness of the realization}
\setcounter{equation}{0}
\label{sec:three}

In this section we assume that the function $K(n,m)=r(n-m)$ has a finite number of negative squares for $n,m\in\mathbb Z$, and prove that
the function $r(n)$ has a representation of the form \eqref{voltaire}. We also prove the uniqueness of this representation under hypothesis
\eqref{chatelet}. We begin with a preliminary proposition.

\begin{Pn}
There exists $C>0$ and $K>0$ such that
\begin{equation}
\label{bound}
\|r(n)\|\le K\cdot C^{|n|},\quad n\in\mathbb Z.
\end{equation}
\end{Pn}

\begin{proof} The claim is true in the scalar
complex-valued case and follows from \eqref{klformula}; see
\cite{MR0107821}. The idea is to reduce the problem to this case.
We write $r(n)=a(n)+jb(n)$ where $a(n)$ and $b(n)$ are $\mathbb
C^{s\times s}$-valued. We obtain a bound of the required form for
every entry of $a(n)$ and $b(n)$. The coefficients $K$ and $C$ in
\eqref{bound} will depend on the given entry. Since there are
$2s^2$ entries, we obtain a bound independent of the
entry.\\

STEP 1: {\sl For every choice of $(e,f)\in\mathbb H^s\times
\mathbb H^s$, the function
\[
K_{e,f}(n,m)=e^*a(n-m)e+f^*\overline{a(n-m)}f+e^*b(n-m)f
-f^*\overline{b(n-m)}e
\]
has at most $2\kappa$ negative squares.}\smallskip

Indeed, the $\mathbb C^{2s\times 2s}$ function
\[
K_1(n,m)=\begin{pmatrix}a(n-m)&b(n-m)\\
-\overline{b(n-m)}&\overline{a(n-m)}\end{pmatrix}
\]
has $2\kappa$ negative squares (See \cite[Proposition 11.4, p.
466]{as3}), and so, for every  fixed choice of $(e,f)\in\mathbb H^s\times \mathbb H^s$, the function
\[
K_{e,f}(n-m)=
\begin{pmatrix}e^*&f^*\end{pmatrix}\begin{pmatrix}a(n-m)&b(n-m)\\
-\overline{b(n-m)}&\overline{a(n-m)}\end{pmatrix}
\begin{pmatrix}e\\f\end{pmatrix}
\]
has at most $2\kappa$ negative squares.\\

STEP 2: {\sl The claim holds for every diagonal entry of
$a(n)$}.\smallskip

Take $e=e_j\in\mathbb H^s$ to be the vector with all entries
equal to $0$, except the $j$-th one equal to $1$ and $f=0$. We
have
\[
K_{e,f}(n,m)=a_{jj}(n-m),
\]
and the result follows from \cite{MR0107821}.\\

STEP 3: {\sl The claim holds for all the entries of
$a(n)$.}\smallskip

Let $\ell\not=j\in\left\{1,\ldots, s\right\}$. We now take
$e=e_{\ell j}(\e)\in\mathbb H^s$ to be the vector with all entries
equal to $0$, except the $\ell$-th one equal to $1$, and the
$j$-th entry equal to $\epsilon$ (where $\varepsilon$ will be
determined) and $f=0$. We have
\[
K_{e,f}(n,m)=a_{\ell \ell
}(n-m)+a_{jj}(n-m)+\overline{\varepsilon}a_{j\ell}(n-m)+\e
a_{\ell j}(n-m).
\]
This function has at most $\kappa$ negative squares and so the
sequence
\[
a_{\ell \ell}(n)+a_{jj}(n)+\overline{\varepsilon}a_{j\ell}(n)+\e
a_{\ell j}(n).
\]
has a bound of the form \eqref{bound} (where $K$ and $C$ depend
on $\ell,j$ and $\epsilon$). The choices $\epsilon=1$ and $\e=i$
gives that the functions
\[
a_{\ell \ell}(n-m)+a_{jj}(n-m)+a_{j\ell}(n-m)+a_{\ell j}(n-m)
\]
and
\[
a_{\ell\ell}(n-m)+a_{jj}(n-m)+i(-a_{j\ell}(n-m)+ a_{\ell j}(n-m))
\]
have at most $\kappa$ negative squares and so the functions
\[
|a_{\ell\ell}(n)+a_{jj}(n)+a_{j\ell}(n)+a_{\ell j}(n)|\le K_1
C_1^{|n|}
\]
and
\[
|a_{\ell\ell}(n)+a_{jj}(n)+i(-a_{j\ell}(n)+ a_{\ell j}(n))|\le K_2
C_2^{|n|}
\]
where the constants depend on $(\ell,j)$. Since $a_{\ell\ell}(n)$
and $a_{jj}(n)$ admit similar bounds we get that both $a_{\ell
j}(n)$
and $a_{j\ell}(n)$ admit bounds of the form \eqref{bound}.\\

STEP 4: {\sl The claim holds for the diagonal entries of
$b(n)$.}\smallskip

We now take $e=e_\ell$ and $f=\e e_j$, where $\e$ is of modulus
$1$. We have
\[
K_{e,f}(n-m)=A(n-m)+B(n-m)
\]
where
\[
\begin{split}
A(n-m)&=a_{\ell\ell}(n-m)+a_{jj}(n-m),\\
B(n-m)&=\overline{\e}b_{\ell\ell}(n-m)+\e b_{jj}(n-m).
\end{split}
\]
The choice $\e=1$ and $\e=i$ lead to the conclusion that
$b_{\ell\ell}(n)$ admits a bound of the form
\eqref{bound} since, as follows from the previous step, $A(n-m)$ admits such a bound.\\

STEP 5: {\sl The claim holds for all the entries of
$b(n)$.}\smallskip

We now take $e=e_{\ell j}(\e_1)$ and $f=\e e_{\ell j}(e_2)$, where
$\e_1$ and $\e_2$ are of modulus $1$. We have now
$K_{e,f}(n-m)=A(n-m)+B(n-m)$  with
\[
\begin{split}
A(n-m)&=e^*a(n-m)e+f^*\overline{a(n-m)}f^*\\
B(n-m)&=e^*b(n-m)f-f^*\overline{b(n-m)}e\\
      &=b_{\ell\ell}(n-m)+\overline{\e_1}\e_2 b_{jj}(n-m)+\overline{\e_1}
      b_{\ell j}(n-m)+\e_2b_{j\ell}(n-m)-\\
&\hspace{5mm}-\overline{b_{\ell\ell}(n-m)}-\overline{\e_2}\e_1
\overline{b_{jj}(n-m)} -\overline{\e_2}\overline{b_{\ell
j}(n-m)}-\e_1\overline{b_{\ell j}(n-m)}.
\end{split}
\]

In view of the previous steps the sequence
\[
\overline{\e_1}b_{j\ell}(n)+\e_2b_{\ell
j}(n)-\overline{\e_2}\overline{b_{j\ell}(n)}-\e_1\overline{b_{\ell
j}(n)}
\]
admits a bound of the form \eqref{bound}. The choices
\[
(\e_1,\e_2)\in\left\{(1,1),(1,-1),(i,i),(i,-i)\right\}
\]
lead to the functions
\[
\begin{split}
&(\overline{b_{j\ell}(n)}-b_{j\ell}(n))+(b_{ \ell j}(n)-\overline{b_{ \ell j}(n)})\\
&(\overline{b_{j\ell}(n)}+b_{j\ell}(n))+(b_{\ell j}(n)+\overline{b_{ \ell j}(n)})\\
&(\overline{b_{j\ell}(n)}-b_{j\ell}(n))+(b_{\ell j}(n)-\overline{b_{\ell j}(n)})\\
&-(\overline{b_{j\ell}(n)}+b_{j\ell}(n))+i(b_{\ell j}(n)-\overline{b_{\ell j}(n)})\\
&-i(\overline{b_{j\ell}(n)}+b_{j\ell}(n))-i(b_{\ell
j}(n)+\overline{b_{\ell j}(n)})
\end{split}
\]
all admit a bound of the form \eqref{bound}.
\end{proof}

\begin{proof}[Proof of Theorem \ref{montmartre}]

 We proceed in a number
of steps to prove the necessity part of the theorem. The first
step is a
direct computation which is omitted.\\

STEP 1: {\sl Let $V$ be a coisometry (that is, $VV^*=I$) in the
quaternionic Pontryagin space $\mathcal P$. Then,
\[
U=\begin{pmatrix}V^*&I-V^*V\\
0&V\end{pmatrix}
\]
is unitary from $\mathcal P^2$ into itself, and is such that
\begin{equation}
V^n=\begin{pmatrix}0&I\end{pmatrix}U^n\begin{pmatrix}0 \\I\end{pmatrix},
\quad n=0,1,2,\ldots
\end{equation}}

STEP 2: {\sl The series
\[
\begin{split}
\varphi(p)&=r(0)+2\sum_{n=1}^\infty p^n r(n),\\
K_\varphi(p,q)&=\sum_{n,m\in\mathbb Z}p^nr(n-m)\overline{q}^m,
\end{split}
\]
converge for $p$ and $q$ in a neighborhood $\Omega$ of the
origin, and it holds that
\begin{equation}
\label{st-ambroise}
K_\varphi(p,q)-pK_\varphi(p,q)\overline{q}=\frac{\varphi(p)+\varphi(q)^*}{2},\quad
p,q\in\Omega.
\end{equation}}\smallskip

The asserted convergences follow from \eqref{bound}, while
\eqref{st-ambroise} is a direct computation.\\

STEP 4: {\sl It holds that
\begin{equation}
\label{opera}
K_\varphi(p,q)=\sum_{n=0}^\infty
p^n\left(\frac{\varphi(p)+\varphi(q)^*}{2}\right)\overline{q}^n
\end{equation}
}\smallskip

This is because equation \eqref{st-ambroise} has a unique solution, and that
the right side of \eqref{opera} solves \eqref{st-ambroise}.\\

STEP 5: {\sl $K_\varphi(p,q)$ is has a finite number of negative
squares in $\Omega$.}\smallskip

Note that for every $N\in\mathbb N$  the function
\[
K_{\varphi,N}(p,q)=\begin{pmatrix}I_s&I_sp&\cdots &
I_sp^N\end{pmatrix}
\mathbb T_N\begin{pmatrix}I_s\\ I_s\overline{q}\\
\vdots\\ I_s\overline{q}^N\end{pmatrix}
\]
has a finite number of negative squares, uniformly bounded by
$\kappa$ in $\Omega$. The claim then follows from
\[
K_\varphi(p,q)=\lim_{N\rightarrow\infty} K_{\varphi,N}(p.q)
\]

STEP 6: {\it There exist a right quaternionic Pontryagin space
$\mathcal P$, a unitary operator  $U\in\mathbf L(\mathcal P)$ and a
linear operator $C\in \mathbf L(\mathbb H^s,\mathcal P)$ such that}
\begin{equation}
\varphi(p)=\frac{CC^*}{2}+\sum_{n=1}^\infty p^{n}C^*U^nC,\quad
p\in\Omega
\end{equation}

Indeed, since the expression in the right side of \eqref{opera}
defines a kernel with a finite number of negative squares, we can
apply \cite[Theorem 10.2]{acs2} to see that there exists a right
quaternionic Pontryagin space $\mathcal P_1$, a coisometric
operator $V\in\mathbf L(\mathcal P_1)$ and a bounded operator
$C_1\in\mathbf L(\mathcal P,\mathcal P_1)$ such that
\[
r(n)=C_1^*V^nC_1,\quad n=0,1,\ldots
\]
We now apply STEP 1 to write
\[
r(n)=C_1^*\begin{pmatrix}0&I\end{pmatrix}U^n\begin{pmatrix}0 \\I\end{pmatrix}C_1,
\quad n=0,1,\ldots,
\]
which concludes the proof with $C=\begin{pmatrix}0 \\ C_1\end{pmatrix}$ and $n\ge 0$.
That the formula still holds for negative $n$ follows from $r(-n)=r(n)^*$ and from the
unitarity of $U$.\\

To conclude the proof, we turn to the uniqueness of the representation \eqref{voltaire}.
Consider two representations \eqref{voltaire},
\[
r(n)=C_1^*U_1^nC_1=C_2^*U_2^nC_2,\quad n\in\mathbb Z,
\]
where $U_1$ and $U_2$ are unitary operators in quaternionic
Pontryagin spaces $\mathcal H_1$ and $\mathcal H_2$ respectively.
Consider the space of pairs
\[
R=\left\{ (U_1^nC_1c\,,\, U_2^nC_2c)\,,\, n\in\mathbb Z,\,\,
c\in\mathbb H^s\right\}.
\]
When condition \eqref{chatelet} is in force for both
representations $R$, defines a linear isometric relation with dense
domain and range, and hence, by the quaternionic version of a theorem of Shmulyan (see \cite[Theorem 7.2]{acs2} and see
\cite[p. 29-30]{adrs} for the complex version of this theorem and for the definition of a linear relation in Pontryagin spaces),
$R$ extends to the graph of a unitary
operator, say $S$, from $\mathcal H_1$ into $\mathcal H_2$:
\[
SU_1^nC_1c=U_2^nC_2c,\quad n\in\mathbb Z,\,\, c\in\mathbb H^s
\]
Setting $n=0$ we get $SC_1=C_2$. Then, taking $n=1$ leads to
$(SU_1)C_1=(U_2S)C_1$, and more generally
\[
(SU_1)U_1^nC_1=(U_2S)U_1^nC_1,\quad n\in\mathbb Z,
\]
and so $SU_1=U_2S$.
\end{proof}

\section{Herglotz's theorem in the quaternionic setting}
\label{Bochner}
\setcounter{equation}{0}
Herglotz's theorem has been already recalled in Section 1, see Theorem \ref{Tm:Oct27yt1}.
Here we state a related result which will be useful in the sequel (see \cite[Theorem 1.3.6]{MR1038803}).
\begin{Tm}
\label{Tm:Oct27wq1}
Let $\mu$ and $\nu$ be $\mathbb C^{s \times s}$-valued measures on $[0, 2\pi ] $. If
$$\int_0^{2\pi} e^{i n t} d\mu(t) = \int_0^{2\pi} e^{i n t} d\nu(t), \quad n \in \mathbb Z,$$
then $\mu = \nu$.
\end{Tm}

Given $P \in \mathbb H^{s \times s}$, there exist unique $P_1, P_2 \in \mathbb C^{s \times s}$ such that $P = P_1 + P_2 j$. Thus there is a bijective homomorphism $\chi : \mathbb H^{s \times s} \to \mathbb C^{2s \times 2s}$ given by
\begin{equation}
\label{eq:Oct27jkl1}
\chi \hspace{0.5mm} P =  \begin{pmatrix} P_1 & P_2 \\ - \overline{P}_2 & \overline{P}_1 \end{pmatrix} \quad {\rm where} \; P = P_1 + P_2 j ,
\end{equation}

\begin{Dn}
\label{def:Oct27j1}
Given an $\mathbb H^{s \times s}$-valued measure $\nu$, write $\nu = \nu_1 + \nu_2 j$, where $\nu_1$ and $\nu_2$ are uniquely determined $\mathbb C^{s \times s}$-valued measures. We call a measure $\nu$ on $[0, 2\pi]$ {\it q-positive} if the $\mathbb C^{2s \times 2s}$-valued measure
\begin{equation}
\label{eq:Nov7kr1}
\mu = \begin{pmatrix} \nu_1 & \nu_2 \\ \nu^*_2 & \nu_3 \end{pmatrix}, \quad {\rm where}\; d\nu_3(t) = d\bar{\nu}_1(2\pi - t),\;\; t \in [0, 2\pi)
\end{equation}
is positive and, in addition,
$$d\nu_2(t) = -d\nu_2(2\pi -t)^T, \quad t \in [0, 2\pi),$$
\end{Dn}

\begin{Rk}
\label{rem:Oct27k1}{\rm
If $\nu$ is $q$-positive, then $\nu = \nu_1 + \nu_2 j$, where $\nu_1$ is a uniquely determined positive $\mathbb C^{s \times s}$-valued measure and $\nu_2$ is a uniquely determined $\mathbb C^{s \times s}$-valued measure.}
\end{Rk}

\begin{Rk}
\label{rem:Nov7uuu1}{\rm
If $r = (r(n))_{n \in \mathbb Z}$ is a $\mathbb H^{s \times s}$-valued sequence on $\mathbb Z$ such that
$$r(n) = \int_0^{2\pi} e^{i n t} d\nu(t),$$
where $\nu$ is a $q$-positive measure, then $r$ is Hermitian, i.e., $r(-n)^* = r(n)$. Indeed,
write $\nu = \nu_1 + \nu_2 j$, where $\nu_1$ and $\nu_2$ are as in Definition \ref{def:Oct27j1}. Then
\begin{align*}
r(-n)^* =& \; \int_0^{2\pi} (d\nu_1(t) - j d\nu_2(t)^*) e^{i n t} \\
=& \; \int_0^{2\pi} e^{i n t} d\nu_1(t) - \int_0^{2\pi} e^{- i n t} (-d\nu_2(t)^T) j \\
=& \; \int_0^{2\pi} e^{i n t}d\nu_1(t) + \int_0^{2\pi} e^{i n t} (-d\nu_2(2\pi - t)^T) j \\
=& \; \int_0^{2\pi} e^{i n t} d\nu_1(t) + \int_0^{2\pi} e^{i n t} d\nu_2(t) \\
=& \; r(n), \quad\quad n \in \mathbb Z
\end{align*}
}
\end{Rk}

\begin{Tm}
\label{Tm:Oct27j1}
The function $n \mapsto r(n)$ from $\mathbb Z$ into $\mathbb H^{s \times s}$ is positive definite if and only if there exists a unique $q$-positive measure $\nu$ on $[0, 2\pi ]$ such that
\begin{equation}
\label{eq:Oct27j1}
r(n) = \int_0^{2\pi} e^{i n t} d\nu(t), \quad n \in \mathbb{Z}.
\end{equation}
\end{Tm}

\begin{proof}
Let $(r(n))_{n \in \mathbb Z}$ be a positive definite sequence and write $r(n) = r_1(n) + r_2(n)j$, where $r_1(n),r_2(n) \in \mathbb C^{s \times s}$, $n \in \mathbb Z$. Put
$R(n) = \chi \hspace{0.5mm} r(n)$, $n \in \mathbb Z$. It is easily seen that $(R(n))_{n \in \mathbb Z}$ is a positive definite $\mathbb C^{2s \times 2s} $-valued sequence if and only if $(r(n))_{n \in \mathbb Z}$ is a positive definite $\mathbb H^{s \times s}$-valued sequence. Thus, by Theorem \ref{Tm:Oct27yt1} there exists a unique positive $\mathbb C^{2s \times 2s}$-valued measure $\mu$ on $[0, 2\pi ] $ such that
	\begin{equation}
\label{eq:Oct27uy1}
R(n) = \int_0^{2\pi} e^{i n t} d\mu(t),\quad n \in \mathbb Z.
\end{equation}
Write $$\mu = \begin{pmatrix} \mu_{11} & \mu_{12} \\ {\mu}^*_{12} & \mu_{22} \end{pmatrix}: \begin{array}{ccc} \mathbb C^{s} & & \mathbb C^s \\ \oplus & \to & \oplus \\ \mathbb C^s & & \mathbb C^s \end{array}.$$
It follows from
$$R(n) = \begin{pmatrix} r_1(n) & r_2(n) \\ -\overline{r_2(n)} & \overline{r_1(n)} \end{pmatrix}, \quad n \in \mathbb Z$$
and \eqref{eq:Oct27uy1} that
$$r_1(n) = \int_0^{2\pi} e^{i n t} d\mu_{11}(t) = \int_0^{2\pi} e^{-i n t} d\bar{\mu}_{22}(t), \quad n \in \mathbb Z$$
and hence
$$\int_0^{2\pi} e^{i n t} d\mu_{11}(t) = \int_0^{2\pi} e^{i n t} d\bar{\mu}_{22}(2\pi - t),\quad n \in \mathbb Z.$$
Thus, Theorem \ref{Tm:Oct27wq1} yields that $d\mu_{11}(t) = d\bar{\mu}_{22}(2\pi -t)$ for $t \in [0, 2\pi)$.
Similarly,
$$r_2(n) = \int_0^{2\pi} e^{i n t } d\mu_{12}(t) = - \int_0^{2\pi} e^{-i n t} d\mu_{12}(t)^T,\quad n \in \mathbb Z$$
and hence
$$\int_0^{2\pi} e^{i n t} d\mu_{12}(t) = \int_0^{2\pi} e^{i n t} (-d\mu_{12}(2\pi - t)^T),\quad n \in \mathbb Z .$$
Thus, Theorem \ref{Tm:Oct27wq1} yields that $d\mu_{12}(t) = -d\mu_{12}(2\pi - t)^T$ for $t \in [0, 2\pi)$.

It is easy to show that
$$\begin{pmatrix} I_s & -j I_s \end{pmatrix} R(n) \begin{pmatrix} I_s \\ j I_s \end{pmatrix} = 2r(n)$$
and hence \eqref{eq:Oct27uy1} yields
\begin{align*}
2r(n) =& \; \int_0^{2\pi} \begin{pmatrix} e^{i n t} & -j e^{i n t} \end{pmatrix} \begin{pmatrix} d\mu_{11}(t) + d\mu_{12}(t) j \\ d\mu_{12}(t)^* + d{\mu}_{22}(t) j \end{pmatrix} \\
=& \; \int_0^{2\pi} e^{i n t} d\mu_{11}(t) + \int_0^{2\pi} e^{i n t} d\mu_{12}(t)j  - \int_0^{2\pi} e^{-i n t} d\mu_{12}(t)^T j \\
& \;\;\;\; + \int_0^{2\pi} e^{-i n t} d\bar{\mu}_{22}(t) \\
=& \; \int_0^{2\pi} e^{i n t} d\mu_{11}(t) + \int_0^{2\pi} e^{i n t} d\mu_{12}(t)j  - \int_0^{2\pi} e^{i n t} d\mu_{12}(2\pi - t)^T j  \\
& \;\;\;\; + \int_0^{2\pi} e^{i n t} d\bar{\mu}_{22}(2\pi - t) \\
=& \; 2 \int_0^{2\pi} e^{i n t} d\mu_{11}(t) + 2 \int_0^{2\pi} e^{i n t} d\mu_{12}(t)j , \quad n \in \mathbb Z,
\end{align*}
where the last line follows from $d\mu_{11}(t) = d\bar{\mu}_{22}(2\pi -t)$ and $d\mu_{12}(t) = -d\mu_{12}(2\pi -t)^T$. If we put $\nu = \mu_{11} + \mu_{12} j$, then $\nu$ is a  $q$-positive measure which satisfies \eqref{eq:Oct27j1}.

Conversely, suppose $\nu = \nu_1 + \nu_2 j$ is a $q$-positive measure on $[0, 2\pi ] $ and put
$$r(n) = \int_0^{2\pi} e^{i n t} d\nu(t),\quad n \in \mathbb Z.$$
Since $\nu$ is $q$-positive,
$$\mu = \begin{pmatrix} \nu_1 & \nu_2 \\ \nu^*_2 & \nu_3 \end{pmatrix},\quad {\rm where}\;\;  d\nu_3(t) = d\bar{\nu}_1(2\pi - t), \;\; t \in [0, 2\pi),$$
is a positive $\mathbb C^{2s \times 2s}$-valued measure on $[0, 2\pi] $ and
$$d\nu_2(t) = -d\nu_2(2\pi - t)^T,\quad t \in [0,2 \pi).$$
Since $\mu$ is a positive $\mathbb C^{2s \times 2s}$-valued measure, $(R(n))_{n \in \mathbb Z}$ is a positive definite $\mathbb C^{2s \times 2s}$-valued sequence, where
$$R(n) := \int_0^{2\pi} e^{i n t } d\mu(t), \quad n \in \mathbb Z,$$
Moreover,
$R(n)$ can be written in form
$$R(n) = \begin{pmatrix} r_1(n) & r_2(n) \\ -\overline{r_2(n)} & \overline{r_1(n)} \end{pmatrix},\quad n \in \mathbb Z,$$
where
\begin{align*}
r_1(n) =& \; \int_0^{2\pi} e^{i n t} d\nu_1(t) ,\quad n \in \mathbb Z{\rm ;} \\
r_2(n) =& \; \int_0^{2\pi} e^{i n t} d\nu_2(t),\quad n \in \mathbb Z.
\end{align*}
Thus, $R(n) = \chi \hspace{0.5mm} r(n)$, where
$$r(n) = r_1(n) + r_2(n) j = \int_0^{2\pi} e^{i n t} d\nu(t).$$
Since $(R(n))_{n \in \mathbb Z}$ is a positive definite $\mathbb C^{2s \times 2s}$-valued sequence we get that $(r(n))_{n \in \mathbb Z}$ is a positive definite $\mathbb H^{s \times s}$-valued sequence.

Finally, suppose that the $q$-positive measure $\nu$ were not unique, i.e., there exists $\tilde{\nu}$ so that $\tilde{\nu} \neq \nu$ and
$$r(n) = \int_0^{2\pi} e^{i n t} d\nu(t) = \int_0^{2\pi} e^{i n t} d\tilde{\nu}(t), \quad n \in \mathbb Z.$$
Write $\nu = \nu_1 + \nu_2 j$ and $ \tilde{\nu} = \tilde{\nu}_1 + \tilde{\nu}_2 j$ as in Remark \ref{rem:Oct27k1}. If we consider $R(n) = \chi \hspace{0.5mm} r(n), n \in \mathbb Z$, then it follows from Theorem \ref{Tm:Oct27yt1} that $\nu_1 = \tilde{\nu}_1$ and $\nu_2 = \tilde{\nu}_2$ and hence that $\nu = \tilde{\nu}$, a contradiction.
\end{proof}

\begin{Rk}{\rm The statement and proof of Herglotz's theorem have been written using an exponential involving the imaginary unit $i$ of the quaternions. Analogous statements can be written using the imaginary units $j$ or $k$ in the basis or with respect to new basis elements chosen in $\mathbb S$.}
\end{Rk}

\section{A theorem of Carath\'eodory in the quaternionic setting}
\label{carath:sec}
\setcounter{equation}{0}

\begin{Dn}
\label{def:Jan4i1}
A function $r: \{-N, \ldots, N \} \to \hh^{s \times s}$ is called {\it positive definite} if $\mathbb{T}_N\succeq 0$, where $\mathbb{T}_N$ is the matrix
defined in (\ref{eq:toeplitz}).
\end{Dn}

\begin{Dn}
Let $r: \{ -N, \ldots, N \} \to \hh^{s \times s}$ be positive definite. We will say that $r$ has a {\it positive definite extension} if there exists a positive definite function $\tilde{r}: \mathbb{Z} \to \hh^{s \times s}$ such that
$$\tilde{r}(n) = r(n), \quad\quad n = -N, \ldots, N.$$
\end{Dn}

\begin{Tm}
\label{thm:Jan4uy1}
If $r: \{ -N, \ldots, N \} \to \hh^{s \times s}$ is positive definite, then $r$ has a positive definite extension.
\end{Tm}

\begin{Rk}{\rm
\label{rem:Jan4t1}
The strategy for proving Theorem \ref{thm:Jan4uy1} is to establish the existence of $r(N+1), r(N+2), \ldots$ so that the block matrices
\begin{align*}
\mathbb{T}_{N+1} =& \; \begin{pmatrix} r(0) & \cdots & r(N+1) \\
\vdots & \ddots & \vdots  \\
r(-N-1) & \cdots & r(0)  \end{pmatrix} \succeq 0 \\
\mathbb{T}_{N+2} =& \; \begin{pmatrix} r(0) & \cdots & r(N+2) \\
\vdots & \ddots & \vdots  \\
r(-N-2) & \cdots & r(0)  \end{pmatrix} \succeq 0, \quad \ldots
\end{align*}
Here we let $r(-N-1) = r(N+1)^*$, $r(-N-2) = r(N+2)^*, \ldots$. We must first establish some lemmas before proving Theorem \ref{thm:Jan4uy1}. The proofs of Lemmas \ref{lem:Jan7u1}, \ref{lem:Jan7r1} and \ref{lem:Jan4qb1} are adapted from Lemma 2.4.2, Corollary 2.4.3 and Theorem 2.4.5 in Bakonyi and Woerdeman \cite{BakonyiWoerdeman}, respectively.}
\end{Rk}

\begin{La}
\label{lem:Jan7u1}
If $A \in \hh^{t \times s}$ and $B \in \hh^{u \times s}$, then
$$B^* B \succeq A^* A$$
if and only if there exists a contraction $G: \ran B \to \ran A$ such that $A = GB$. Moreover,
$G$ is unique and an isometry if and only if $B^* B = A^* A$.
\end{La}

\begin{proof}
If there exists a contraction $G: \ran B \to \ran A$ such that $A = GB$, then it is easy to verify that $B^* B \succeq A^* A$. Conversely, if $B^* B \succeq A^* A$, then let $y \in \ran B$, i.e $y = Bx$ for some $x \in \hh^{s}$. Let $G: \ran B \to \ran A$ be given by
$$Gy = Ax.$$
To check that $G$ is well-defined, suppose that
$$y = Bx = B \tilde{x},$$
where $\tilde{x} \in \hh^s$. Using $B^* B \succeq A^* A$ we get that
$$0 \leq (x - \tilde{x})^* A^* A (x- \tilde{x}) \leq (x - \tilde{x})^* B^* B (x - \tilde{x}) = 0$$
and hence $Ax = A \tilde{x}$. Therefore, $G$ is well-defined.

We will now show that $G$ is a contraction. Let $\{ y_n \}_{n=1}^{\infty}$ be a convergent sequence in $\ran B$. If $\{ x_n \}_{n=1}^{\infty}$ in $\hh^s$ so that
$$B x_n = y_n,$$
then
\begin{align}
(G y_n - G y_m )^* (G y_n - G y_m ) =& \; [A(x_n - x_m)]^*  [A(x_n - x_m)] \nonumber \\
\leq& \; [B(x_n - x_m)]^*  [B(x_n - x_m)] \nonumber \\
=& \; (y_n - y_m)^* (y_n - y_m). \label{eq:Jan7t1}
\end{align}
Since $\{ y_n \}_{n=1}^{\infty}$ is a convergent sequence in $\ran B$, $\{ y_n \}_{n=1}^{\infty}$ is also a Cauchy sequence in $\ran B$ and hence $\{ G y_n \}_{n=1}^{\infty}$ is a Cauchy sequence as well. Thus,
$$\lim_{n\uparrow \infty} G y_n$$
exists. The inequality given in \eqref{eq:Jan7t1} readily yields that $y^* G^* G y \leq y^* y$, whence $G$ is a contraction. Note that $G$ is unique by construction, since
the equation $A = GB$ requires that whenever $y = B x $ we get that $G y = Ax$.

Finally, if $G$ is an isometry then it follows from the equality $A = GB$ that $A^* A = B^* B$. Conversely, if $A^* A = B^* B$, then $ y = Bx$ and $G y = Ax$ yield that
$$ y^* G^* G y = x^* A^* A x = x^* B^* B x = y^* y. $$
Thus, $u^* G^* G y = y^* y$ for $y \in \ran B$.
\end{proof}

\begin{La}
\label{lem:Jan7r1}
If
$$K = \begin{pmatrix} A & B \\ B^* & C \end{pmatrix} \in \hh^{ (t+u) \times (s+u) },$$
then $K \succeq 0$ if and only if the following conditions hold{\rm :}
\begin{enumerate}
\item[(i)] $A \succeq 0$ and $C \succeq 0${\rm ;}
\item[(ii)] $B = A^{1/2} G C^{1/2}$ for some contraction $G: \ran C \to \ran A$.
\end{enumerate}
\end{La}

\begin{proof}
Suppose conditions (i) and (ii) are in force. It follows from (i) that there exist $P$ and $Q$ such that $A = P^* P$ and $C = Q^* Q$. Thus,
$$K = \begin{pmatrix} P^* & 0 \\ 0 & Q^* \end{pmatrix} \begin{pmatrix} I & G \\ G^* & I \end{pmatrix} \begin{pmatrix} P & 0 \\ 0 & Q \end{pmatrix} \succeq 0,$$
since $G$ is a contraction. Conversely, suppose $K \succeq 0$ and let $P$ and $Q$ be given by
$$ \begin{pmatrix} P^* \\ Q^* \end{pmatrix} \begin{pmatrix} P & Q \end{pmatrix} = \begin{pmatrix} A & B \\ B^* & C \end{pmatrix}.$$
Thus, $P^* P = A^{1/2} A^{1/2}$ and $ Q^* Q = C^{1/2} C^{1/2}$. Using Lemma \ref{lem:Jan7u1} we arrive at the isometries $G_1: \ran A \to \ran P$ and $G_2: \ran C \to \ran Q$ which satisfy
$P = G_1 A^{1/2}$ and $Q = G_2 C^{1/2}$, respectively. Therefore,
$$B = P^* Q = A^{1/2} G_1^* G_2 C^{1/2}$$
and thus $B = A^{1/2} G C^{1/2}$, where $G = G_1^* G_2$ is a contraction.
\end{proof}

\begin{Dn}
\label{def:Jan4tq1}
We will call a block matrix, with quaternionic entries,
$$K = \begin{pmatrix} A & B & ? \\ B^* & C  & D \\ ? & D^* & E \end{pmatrix}$$
{\it partially positive semidefinite} if all principle specified minors are nonnegative. We will say that $K$ has a {\it positive semidefinite completion} if there exists a quaternionic matrix $X$ so that
$$\begin{pmatrix} A & B & X \\ B^* & C & D \\ X^* & D^* & E \end{pmatrix} \succeq 0.$$
\end{Dn}

\begin{La}
\label{lem:Jan4qb1}
If
$$K = \begin{pmatrix} A & B & ? \\ B^* & C & D \\ ? & D^* & E \end{pmatrix}$$
is partially positive semidefinite, then $K$ has a positive semidefinite completion given as follows. Let $G_1: \ran C \to \ran A$ and $G_2: \ran E \to \ran C$ be contractions so that $B = A^{1/2} G_1 C^{1/2}$and $D = C^{1/2} G_2 E^{1/2}$. Choosing the $(1,3)$ block entry of $K$ to be $A^{1/2} G_1 G_2 E^{1/2}$ results in a positive semidefinite completion.
\end{La}

\begin{proof}
Since $K$ is partially positive semidefinite,
$$K_1  = \begin{pmatrix} A & B \\ B^* & C \end{pmatrix} \succeq 0 \quad\quad
{\rm and} \quad\quad K_2 = \begin{pmatrix} C & D \\ D^* & E \end{pmatrix} \succeq 0.$$
Use Lemma \ref{lem:Jan7r1} on $K_1$ and $K_2$ to produce contractions $G_1$ and $G_2$, resepectively, so that $B = A^{1/2} G_1 C^{1/2}$ and $D = C^{1/2} G_2 E^{1/2}$. Since $G_1$ and $G_2$ are contractions, the factorization
\begin{align*}
& \;\;\;\;\;\;\;\;\;\;\;\;\;\;\;\; \widetilde{K} = \begin{pmatrix} A & B & A^{1/2} G_1 G_2 E^{1/2} \\
B^* & C & D \\
E^{1/2} (G_2)^* (G_1)^* A^{1/2} & D^* & E \end{pmatrix} \\
=& \; \begin{pmatrix} A & A^{1/2} G_1 C^{1/2} & A^{1/2} G_1 G_2 E^{1/2} \\
C^{1/2} (G_1)^* A^{1/2} & C & C^{1/2} G_2 E^{1/2} \\
E^{1/2} (G_2)^* (G_1)^* A^{1/2} & E^{1/2} (G_2)^* C^{1/2} & E \end{pmatrix} \\
=& \; \begin{pmatrix} A^{1/2} & 0 & 0 \\ 0 & C^{1/2} & 0 \\ 0 & 0 & E^{1/2} \end{pmatrix}
\begin{pmatrix} I & G_1 & G_1 G_2 \\ (G_1)^* & I & G_2 \\ (G_2)^*(G_1)^* & (G_2)^* & I \end{pmatrix} \begin{pmatrix} A^{1/2} & 0 & 0 \\ 0 & C^{1/2} & 0 \\ 0 & 0 & E^{1/2} \end{pmatrix} \\
=& \; \begin{pmatrix} A^{1/2} & 0 & 0 \\ 0 & C^{1/2} & 0 \\ 0 & 0 & E^{1/2} \end{pmatrix}
\begin{pmatrix} I & 0 & 0 \\ (G_1)^* & I & 0 \\ (G_1G_2)^* & (G_2)^* & I \end{pmatrix}
\begin{pmatrix} I & 0 & 0 \\ 0 & I - (G_1)^*G_1 & 0 \\ 0 & 0 & I - (G_2)^* G_2 \end{pmatrix} \\
& \;\;\;\;\;\;\;\;\;\;\; \times \begin{pmatrix} I & G_1 & G_1 G_2 \\ 0 & I & G_2 \\ 0 & 0 & I \end{pmatrix} \begin{pmatrix} A^{1/2} & 0 & 0 \\ 0 & C^{1/2} & 0 \\ 0 & 0 & E^{1/2} \end{pmatrix},
\end{align*}
shows that $\widetilde{K}$ is a positive semidefinite completion of $K$.
\end{proof}

We are now ready to prove Theorem \ref{thm:Jan4uy1}.

\begin{proof}[Proof of Theorem \ref{thm:Jan4uy1}]
Let $r: \{ -N, \ldots, N \} \to \hh^{s \times s}$ be positive definite. It follows from Lemma \ref{lem:Jan4qb1} with
\begin{align*}
A =& \; r(0){\rm ;} \\
B =& \; \begin{pmatrix} r(1) & \cdots & r(N) \end{pmatrix}{\rm ;} \\
C =& \; \begin{pmatrix} r(0) & \cdots & r(N) \\ \vdots & \ddots & \vdots \\ r(-N) & \cdots & r(0) \end{pmatrix} {\rm ;} \\
D =& \; \begin{pmatrix} r(N)^T & \cdots & r(1)^T \end{pmatrix}^T{\rm ;} \\
E =& \; A,
\end{align*}
that there exist contractions $G_1$ and $G_2$ so that if we put $r(N+1) = A^{1/2} G_1 G_2 E^{1/2}$ and $r(-N-1) = r(N+1)^*$, then
$$\begin{pmatrix}
r(0) & \cdots & r(N+1)  \\
\vdots & \ddots & \vdots  \\
r(-N-1) & \cdots & r(0)
\end{pmatrix} \succeq 0.
$$
Continuing in this fashion, we can choose $r(N+2), r(N+3), \ldots$ so that
$$\begin{pmatrix} r(0) & \cdots & r(N+2) \\ \vdots & \ddots & \vdots \\ r(-N-2) & \cdots & r(0)\end{pmatrix} \succeq 0, \quad \begin{pmatrix} r(0) & \cdots & r(N+3) \\ \vdots & \ddots & \vdots \\ r(-N-3) & \cdots & r(0)\end{pmatrix} \succeq 0, \ldots.$$
Thus we have contructed $\tilde{r}: \mathbb{Z} \to \hh^{s \times s}$ which is positive definite and satisfies
$$\tilde{r}(n) = r(n), \quad\quad n = -N, \ldots, N.$$
\end{proof}

\section{A theorem of Krein and Iohvidov in the quaternionic setting}\label{IohvidovKrein}
Sasv\'ari \cite{sasvari} attributes the following theorem to Krein and Iohvidov \cite{MR0107821}.

\begin{Tm}
\label{Tm:Nov7k1}
Let $a = (a(n))_{n \in \mathbb Z}$ be a bounded Hermitian complex-valued sequence on $\mathbb Z$. The sequence $a$ has $\kappa$ negative squares if and only if
there exist measures $\mu_+$ and $\mu_-$ on $[0, 2\pi]$ and mutually distinct points $t_1, \ldots, t_\kappa \in [0, 2\pi]$ satisfying
$\mu_+(t_j) = 0$ for $j =1, \ldots, \kappa$ and ${\rm supp}\hspace{0.5mm}(\mu_-) = \{ t_1, \ldots, t_k \}$ and such that
\begin{equation}
\label{eq:Nov7a1}
a(n) = \int_0^{2\pi} e^{i n t} d\mu_+(t) - \int_0^{2\pi} e^{i n t} d\mu_-(t),\quad\quad n \in \mathbb Z.
\end{equation}
\end{Tm}

\begin{proof}
A proof for this result when $\mathbb Z$ is replaced by an arbitrary locally compact Abelian group can be found in \cite{sasvari}.
\end{proof}

It will be our goal in this section to obtain a direct analogue of Theorem \ref{Tm:Nov7k1} when $(a(n))_{n \in \mathbb Z}$ is a bounded Hermitian $\mathbb H^{s \times s}$-valued sequence. To achieve this goal, we will first generalize Theorem \ref{Tm:Nov7k1} to the case when $(a(n))_{n \in \mathbb Z}$ is $\mathbb C^{s \times s}$-valued sequence and then the desired result will follow.

\begin{Dn}
\label{Dn:Nov7k1}
Let $ M = \sum_{q=1}^k P_q \delta_{t_q}$ be a $\mathbb C^{s \times s}$-valued measure on $[0, 2\pi ] $, where $\delta_t$ denotes the usual Dirac point measure at $t$. We let
$${\rm card}\hspace{0.5mm}{\rm supp}\hspace{0.5mm}  M = \sum_{q=1}^k {\rm rank}\hspace{0.5mm} P_q.$$
\end{Dn}

\begin{Tm}
\label{Tm:Nov7i1}
Let $A = (A(n))_{n \in \mathbb Z}$ be a bounded Hermitian $\mathbb C^{s \times s}$-valued sequence on $\mathbb Z$. $A$ has $\kappa$ negative squares if and only if there exist positive $\mathbb C^{s \times s}$-valued measures $ M_+$ and $ M_-$ on $[0, 2\pi]$ and mutually distinct points $t_1, \ldots, t_k \in [0, 2\pi ] $
satisfying
$ M_+(t_j) = 0$ for $j =1, \ldots, k$, ${\rm supp}\hspace{0.5mm}( M_-) = \{ t_1, \ldots, t_k \}$ and ${\rm card}\hspace{0.5mm}{\rm supp}\hspace{0.5mm}  M_- = \kappa$ and such that
\begin{equation}
\label{eq:Nov7a355}
A(n) = \int_0^{2\pi} e^{i n t} d M_+(t) - \int_0^{2\pi} e^{i n t} d M_-(t),\quad\quad n \in \mathbb Z.
\end{equation}
\end{Tm}

\begin{proof}
If $A$ has $\kappa$ negative squares, then $a_v = (v^* A(n) v )_{n \in \mathbb Z}$ will be a complex-valued sequence with at most $\kappa$ negative squares for any $v \in \mathbb C^s$. It follows then from Theorem \ref{Tm:Nov7k1} that there exist measures $\mu_+^{(v)}$ and $\mu_-^{(v)}$ on $[0, 2\pi] $ and mutually distinct points $t_1^{(v)}, \ldots, t_{k_v}^{(v)} \in [0, 2\pi]$ satisfying
$\mu_+^{(v)}(t_j^{(v)}) = 0$ for $j =1, \ldots, k_v$ and ${\rm supp}\hspace{0.5mm}(\mu_-^{(v)}) = \{ t_1^{(v)}, \ldots, t_{k_v}^{(v)} \}$, where $k_v \leq \kappa$, and such that
\begin{equation}
a_v(n) = \int_0^{2\pi} e^{i n t} d\mu_+^{(v)}(t) - \int_0^{2\pi} e^{i n t} d\mu_-^{(v)}(t),\quad\quad n \in \mathbb Z.
\end{equation}
Let
$$4 \mu_{\pm}^{(v,w)} = \mu_{\pm}^{(v+w)} - \mu_{\pm}^{(v-w)} + i \mu_{\pm}^{(v + i w)} - i \mu_{\pm}^{(v-iw)},\quad\quad v,w \in \mathbb C^s.$$
Then there exist positive $\mathbb C^{s \times s}$-valued measures $ M_{\pm}$ such that
$\langle  M_{\pm} v, w \rangle = \mu_{\pm}^{(v,w)}$ and
\begin{equation}
\label{eq:Oct7ju1}
A(n) = \int_0^{2\pi} e^{i n t} dM_+(t) - \int_0^{2\pi} e^{i n t} d M_-(t), \quad\quad n \in \mathbb Z.
\end{equation}
It follows from \eqref{eq:Oct7ju1} together with the fact that $A$ has $\kappa$ negative squares that $${\rm card}\hspace{0.5mm} {\rm supp}\hspace{0.5mm} M_- = \kappa.$$ By construction, $M_+(t_j) = 0$ for all $t_j \in {\rm supp}\hspace{0.5mm} M_+$.

Conversely, suppose \eqref{eq:Nov7a355} is in force. It is easy to check that $A$ is a bounded Hermitian sequence with at most $\kappa$ negative squares. The fact that $A$ has exactly $\kappa$ negative follows from the uniqueness of the measure $ M = M_+ - M_-$ in \eqref{eq:Nov7a3} (see Theorem \ref{Tm:Oct27wq1}).
\end{proof}

\begin{Dn}
\label{Dn:Nov7k7}
Let $ \nu = \nu_1 + \nu_2 j$ be a $q$-positive measure on $[0, 2\pi ] $ with finite support. We let
$${\rm card}\hspace{0.5mm}{\rm supp}\hspace{0.5mm}  \nu = (1/2){\rm card}\hspace{0.5mm}{\rm supp} \hspace{0.5mm} \mu,$$
where $\mu$ is an in \eqref{eq:Nov7kr1}.
\end{Dn}

\begin{Tm}
\label{Tm:Nov7i4}
Let $a = (a(n))_{n \in \mathbb Z}$ be a bounded Hermitian $\mathbb H^{s \times s}$-valued sequence on $\mathbb Z$. The sequence $a$ has $\kappa$ negative squares if and only if there exist $q$-positive measures $ \nu_+$ and $ \nu_-$ on $[0, 2\pi ] $ and mutually distinct points $t_1, \ldots, t_k \in [0, 2\pi ] $ satisfying
satisfying
$ \nu_+(t_j) = 0$ for $j =1, \ldots, k$, ${\rm supp}\hspace{0.5mm}(d \nu_-) = \{ t_1, \ldots, t_k \}$ and ${\rm card}\hspace{0.5mm}{\rm supp}\hspace{0.5mm}  \nu_- = \kappa$ and such that
\begin{equation}
\label{eq:Nov7a8}
a(n) = \int_0^{2\pi} e^{i n t} d \nu_+(t) - \int_0^{2\pi} e^{i n t} d \nu_-(t),\quad\quad n \in \mathbb Z.
\end{equation}
\end{Tm}

\begin{proof}
If $a$ is a bounded Hermitian $\mathbb H^{s \times s}$-valued sequence with $\kappa$ negative squares, then $A = (A(n))_{n \in \mathbb Z}$, where $A(n) = \chi \hspace{0.5mm} a(n)$, has $2\kappa$ negative squares (see Proposition 11.4 in \cite{as3}).  Thus, Theorem \ref{Tm:Nov7i1} guarantees the existence of positive $\mathbb C^{2s \times 2s}$-valued measures $ M_+$ and $ M_-$ on $[0, 2\pi ] $ and mutually distinct points $t_1, \ldots, t_k \in [0, 2\pi ] $ satisfying
satisfying
$ M_+(t_j) = 0$ for $j =1, \ldots, k$, ${\rm supp}\hspace{0.5mm}( M_-) = \{ t_1, \ldots, t_k \}$ and ${\rm card}\hspace{0.5mm}{\rm supp}\hspace{0.5mm} M_- = 2\kappa$ and such that
\begin{equation}
\label{eq:Nov7a3}
A(n) = \int_0^{2\pi} e^{i n t} d M_+(t) - \int_0^{2\pi} e^{i n t} d M_-(t),\quad\quad n \in \mathbb Z.
\end{equation}
If we write
$$dM_{\pm} = \begin{pmatrix} dM_{\pm}^{(11)} & dM_{\pm}^{(12)} \\ (dM_{\pm}^{(12)})^* & dM_{\pm}^{(22)} \end{pmatrix} \begin{array}{ccc} \mathbb C^s & & \mathbb C^s \\ \oplus & \to & \oplus \\ \mathbb C^s & & \mathbb C^s \end{array}$$
and proceed as in the proof of Theorem \ref{Tm:Oct27j1} we get that
$$dM_+^{(11)}(t) - dM_-^{(11)}(t) = dM_+^{(22)}(2\pi -t) - dM_-^{(22)}(2\pi - t),\quad\quad t \in [0, 2\pi)$$
and
$$dM_+^{(12)}(t) - dM_-^{(12)}(t) = -(dM_+^{(12)}(2\pi -t)^T - dM_-^{(12)}(2\pi - t)^T),\quad\quad t \in [0, 2\pi).$$
Consequently, it follows from $d M_+(t_j) = 0$ for $j =1, \ldots, k$ and ${\rm supp}\hspace{0.5mm}(d M_-) = \{ t_1, \ldots, t_k \}$ that
$$dM_-^{(11)}(t) = dM_-^{(22)}(2\pi - t),\quad\quad t \in [0, 2\pi)$$
and
$$dM_-^{(12)}(t) = -dM_-^{(12)}(2\pi - t)^T,\quad\quad t \in [0, 2\pi).$$
Thus,
$$dM_+^{(11)}(t) = dM_+^{(22)}(2\pi - t),\quad\quad t \in [0, 2\pi)$$
and
$$dM_+^{(12)}(t) = -dM_+^{(12)}(2\pi - t)^T,\quad\quad t \in [0, 2\pi).$$
Taking advantage of the above equalities we can obtain
$$a(n) = \int_0^{2\pi} e^{i n t} d\nu_+(t) - \int_0^{2\pi} e^{i n t} d\nu_-(t),\quad \quad n \in \mathbb Z,$$
where $\nu_{\pm}(t) = M_{\pm}^{(11)}(t) + M_{\pm}^{(12)}(t) j$. It is readily checked that $\nu_{\pm}$ are $q$-positive measures. Moreover,
$a$ has $\kappa$ negative squares since $A$ has $2\kappa$ negative squares and $\nu_+(t_j) = 0$ for all $t_j \in {\rm supp}\hspace{0.5mm}\nu_-$ and ${\rm card}\hspace{0.5mm}{\rm supp}\hspace{0.5mm} \nu_- = \kappa$.

Conversely, suppose that $a$ is a $\mathbb H^{s \times s}$-valued sequence which obeys \eqref{eq:Nov7a8}. Consequently, $a$ is bounded. The fact that $a$ is Hermitian follows from Remark \ref{rem:Nov7uuu1}. To see that $a$ has $\kappa$ negative squares, one can consider the $\mathbb C^{2s \times 2s}$-valued sequence $(A(n))_{n \in \mathbb Z}$, where $A(n) = \chi \hspace{0.5mm} a(n)$, $n \in \mathbb Z$ and use the converse statement in Theorem \ref{Tm:Nov7i1} to see that $A$ has $2\kappa$ negative squares. The fact that $a$ has $\kappa$ negative squares then follows by definition.
\end{proof}

\section{Herglotz's integral representation theorem in the scalar case}
\setcounter{equation}{0}
\label{sec:six}
In this section we present an analogue of Herglotz's theorem in the quaternionic scalar case. Even though this is a byproduct of the preceding discussion, it may be useful to have the result stated for scalar valued slice hyperholomorphic functions. We begin by proving an integral representation formula which holds on
$\mathbb{B}_r=\{p\in \mathbb{H}\ : \ |p|<r\}$,
namely on the quaternionic open ball centered at $0$ and with radius $r>0$.
A similar formula which is based on a different representation of a
slice hyperholomorphic function, less useful to determine the real part of a function, is discussed in \cite{CGS_INT_REP}.
Note also that, unlike what happens in the complex case, the real part of a slice hyperholomorphic function is not harmonic.
\begin{La}\label{TEOR1}
Let $f:\mathbb{B}_{1+\varepsilon}\to \mathbb{H}$ be a slice hyperholomorphic function, for some $\varepsilon >0$. Let $I,J\in\mathbb S$ with $J$ orthogonal to $I$ and let $F,G:  \mathbb{B}_{1+\varepsilon}\cap \mathbb{C}_I \to \mathbb{C}_I$ be holomorphic functions
 such that for any $z=x+Iy$ the restriction $f_I$ can be written as
$f_I(z)=F(z)+G(z)J.$
Then, on $\mathbb B_1\cap\mathbb C_I$ the following formula holds:
\[
f_I(z)
=I [{\rm Im} F(0) +{\rm Im} G(0)J]
+\frac{1}{2\pi}\int_0^{2\pi}\frac{e^{It}+z}{e^{It}-z} \, [{\rm Re}(F(e^{It}))+{\rm Re}(G(e^{It}))J] dt.
\]
Moreover
\begin{equation}\label{REPOSITIVE}
{\rm Re }\Big(\frac{e^{It}+z}{e^{It}-z} \, [{\rm Re}(F(e^{It}))+{\rm Re}(G(e^{It}))J]\Big)=
\frac{1-|z|^2}{|e^{It}-z|^2}{\rm Re}(F(e^{It})).
\end{equation}
\end{La}
\begin{proof}
The proof is an easy consequence of the Splitting Lemma \ref{splitting}: for
every fixed $I$, $J \in \mathbb{S}$ such that $J$ is orthogonal to $I$, there are two holomorphic functions
$F,G:  \mathbb{B}_{1+\varepsilon}\cap \mathbb{C}_I \to \mathbb{C}_I$ such that for any $z=x+Iy$, it is
$f_I(z)=F(z)+G(z)J.$
It is immediate that these two holomorphic functions $F$, $G$ satisfy (see p. 206 in \cite{CAPB})
$$
F(z)= I{\rm Im}\, F(0)+\frac{1}{2\pi}\int_0^{2\pi}\frac{e^{It}+z}{e^{It}-z} \, {\rm Re}(F(e^{It})) dt \ \ \ z\in \mathbb{B}_{1}\cap \mathbb{C}_I,
$$
$$
G(z)= I{\rm Im}\,G(0)+\frac{1}{2\pi}\int_0^{2\pi}\frac{e^{It}+z}{e^{It}-z} \, {\rm Re}(G(e^{It})) dt \ \ \ z\in \mathbb{B}_{1}\cap \mathbb{C}_I,
$$
and the first part of the statement follows.
The second part is a consequence of the equality
$$
\frac{e^{It}+z}{e^{It}-z}=\frac{1-|z|^2}{|e^{It}-z|^2}+2I\frac{y\cos t-x\sin t}{|e^{It}-z|^2}
$$
leading to
$$
{\rm Re }\Big(\frac{e^{It}+z}{e^{It}-z} \, [{\rm Re}(F(e^{It}))+{\rm Re}(G(e^{It}))J]\Big)=
\frac{1-|z|^2}{|e^{It}-z|^2}{\rm Re}(F(e^{It})).
$$
\end{proof}

\begin{Rk}\label{realf} {\rm Let $f: \Omega\to\mathbb H$ be a slice hyperholomorphic function and write
 $$
 f(p)=f_0(x_0,\ldots, x_3)+f_1(x_0,\ldots, x_3)i+f_2(x_0,\ldots, x_3)j+f_3(x_0,\ldots, x_3)k,
  $$
  with $f_\ell:\Omega\to\mathbb R$, $\ell=0,\ldots,3$, $p=x_0+x_1i+x_2j+x_3k$. It is easily seen that the restriction $f_i=f_{|\mathbb C_i}$  can be written as
\[
\begin{split}
f_i(x+iy)&=(f_0(x+iy)+f_1(x+iy)i)+(f_2(x+iy)+f_3(x+iy)i)j
\\
&
=F(x+iy)+G(x+iy)j
\end{split}
\]
 and so
$${\rm Re}(f_{|\mathbb C_i})(x+iy)={f_0}_{|\mathbb C_i}(x+iy))={\rm Re}(F)(x+iy).$$
More in general, consider $I,J\in\mathbb S$ with $I$ orthogonal to $J$, and rewrite $i,j,k$ in terms of the imaginary units $I,J,IJ=K$. Then
$$
f(p)=f_0(x_0,\ldots, x'_3)+\tilde f_1(x_0,\ldots, x'_3)I+\tilde f_2(x_0,\ldots, x'_3)J+\tilde f_3(x_0,\ldots, x'_3)K,
 $$
 where $p=x_0+x'_1I+x'_2J+ x'_3K$ and the $x'_\ell$ are linear combinations of the $x_\ell$, $\ell=1,2,3$. The restriction of $f$ to the complex plane $\mathbb C_I$ is then $f_I(x+Iy)=\tilde F(x+Iy)+\tilde G(x+Iy)J$ and reasoning as above we have
$$
{\rm Re}(f_I(x+Iy))={f_0}_{|\mathbb C_I}(x_0,\ldots, x_3')={\rm Re}(\tilde F(x+Iy)).
$$
We conclude that the real part of the restriction $f_I$ of $f$ to a complex plane $\mathbb C_I$ is the restriction of $f_0$ to the given complex plane. Thus if ${\rm Re}(f)$ is positive also the real part of the restriction $f_I$ to any complex plane is positive.
}
\end{Rk}

\begin{Tm}[Herglotz's theorem on a slice]\label{Herglotz}
Let $f:\mathbb{B}_{1}\to \mathbb{H}$ be a slice hyperholomorphic function with ${\rm Re}(f(p))\geq 0$ in $\mathbb B_1$.
Fix $I,J\in\mathbb S$ with $J$ be orthogonal to $I$.
Let $f_I$ be the restriction of $f$ to the complex plane $\mathbb{C}_I$ and let
$F,G:  \mathbb{B}_{1}\cap \mathbb{C}_I \to \mathbb{C}_I$ be holomorphic functions
 such that for any $z=x+Iy$, it is
$f_I(z)=F(z)+G(z)J.$
Then $f_I$ can be written in $\mathbb B_1\cap\mathbb C_I$ as
\begin{equation}\label{Herglotzslice}
f_I(z)=I [{\rm Im} F(0) +{\rm Im} G(0)J]+\int_0^{2\pi}\frac{e^{It}+z}{e^{It}-z} \, d\mu_J(t),
\end{equation}
where $\mu_J(t)=\mu_1+\mu_2 J$ is a finite variation complex measure on $\mathbb C_J$ with $\mu_1$ positive and $\mu_2$ real and of finite variation on $[0,2\pi]$.
\end{Tm}
\begin{proof} The proof follows \cite[p. 207]{CAPB}. First, we note that by Remark \ref{realf}, ${\rm Re}(f(p))\geq 0$ in $\mathbb B_1$ implies that ${\rm Re}(f_I(z))\geq 0$ for $z\in\mathbb B_1\cap\mathbb C_I$.
Let $\rho$ be a real number such that $0<\rho <1$. Then $f_I(\rho z)$ is slice hyperholomorphic in the disc $|z|<1/\rho$ and so by
Lemma \ref{TEOR1} the restriction $f_I(z)$ may be written in $|z|<1$ as
 $$
 f_I(\rho z)
=I [{\rm Im} F(0) +{\rm Im} G(0)J]
+\frac{1}{2\pi}\int_0^{2\pi}\frac{e^{It}+z}{e^{It}-z} \, [{\rm Re}(F(\rho e^{It}))+{\rm Re}(G(\rho e^{It}))J] dt.
 $$
where
 $$
d\mu_J(t,\rho )= \frac{1}{2\pi}[{\rm Re}(F(\rho e^{It}))+{\rm Re}(G(\rho e^{It}))J] dt
 $$
has real positive part, since it is immediate that ${\rm Re}(f_I(\rho e^{It}))={\rm Re}(F(\rho e^{It}))$ and
 $$
 \int_0^{2\pi}d\mu_{J}(t,\rho )= [{\rm Re} F(0) +{\rm Re} G(0)J].
 $$
Let us set
 $$
\Lambda_I(z,t):= \frac{e^{It}+z}{e^{It}-z}
 $$
 and consider
 $$
 \int_0^{2\pi}\Lambda_I(z,t)d\mu_J(t;\rho).
 $$
 Let $\{\rho_n\}_{n\in \mathbb{N}}$ be  a sequence of real numbers  with $0<\rho_n<1$ such that $\rho_n\to 1$ when $n$ goes to infinity.
 To conclude the proof we need Helly's theorem in the complex case. This result assures that the family of finite variation real-valued $d\nu(t;\rho_n)$
 contains a convergent subsequence which tends to $d\nu(t)$ which is of finite variation, in the sense that
 $$
\lim_{n\to \infty} \int_0^{2\pi}\lambda(w,t) d\nu(t,\rho_n)= \int_0^{2\pi}\lambda(w,t) d\nu(t)
 $$
 for every continuous complex-valued function $\lambda(w,t)$.
 In the slice hyperholomorphic setting the integrand is the product of the continuous $\mathbb{C}_I$-valued function $\Lambda_I(z,t)=\Lambda_1(z,t)+I\Lambda_2(z,t)$ where $\Lambda_1$ and $\Lambda_2$ are real valued, and of the $\mathbb{C}_J$-valued
  $d\mu_{J}(t,\rho_n)=  d\mu_1(t,\rho_n)+d\mu_2(t,\rho_n) J$ (since both $d\mu_1(t,\rho_n)$ and  $d\mu_2(t,\rho_n)$ are real-valued).\\
  Then $\Lambda_I(z,t)d\mu_{J}(t)$ can be split in components to which we apply Helly's theorem. The positivity of $d\mu_1$ follows from the positivity of $d\mu_1(t,\rho_n)$, and this completes the proof.

\end{proof}

\begin{Cy}
Let $f$ be slice hyperholomorphic function  on $\mathbb{B}_1$  such that $f(0)=1$.
Suppose that $f$ has real positive part on $\mathbb{B}_1$.
Then its restriction $f_I$ can be represented as
$$
f_I(z)=\int_0^{2\pi}\frac{e^{It}+z}{e^{It}-z} \, d\mu_J(t),
$$
where $\mu_J(t)=\mu_1(t)+\mu_2(t)J$ is a finite variation complex measure on $\mathbb C_J$ with $\mu_1(t)$ positive for $t\in [0,2\pi ] $.
Moreover, the power series expansion of $f$
$$
f(p)=1+\sum_{n=1}^\infty p^n a_n
$$
is such that $|a_n|\leq k$, for some $k\in\mathbb R$, for every $n\in \mathbb{N}$.
\end{Cy}
\begin{proof}
The first part of the corollary immediately follows from Theorem \ref{Herglotz}.
Then observe that
$$
\frac{e^{It}+z}{e^{It}-z}=1+2\sum_{n=1}^\infty z^n e^{-Int}
$$
and so the coefficients $a_n$ in the power series expansion are given by
\begin{equation}\label{an}
a_n=2\int_0^{2\pi} e^{-Int}d\mu_J(t).
\end{equation}
Moreover
$$
|a_n|\leq 2 \int_0^{2\pi}|d\mu_{J}(t)|\leq k,
 $$
 for some $k\in\mathbb R$ since $d\mu_{J}(t)$ is of finite variation and so it is bounded.
\end{proof}
\begin{Rk}{\rm Formula (\ref{an}) expresses $a_n$ in integral form. However, there is an infinite number of ways of writing $a_n$ with a similar expression, depending on the choices of $I$ and $J$ made to write (\ref{Herglotzslice}). An important difference with the result in Section \ref{Bochner} is that the measure $d\mu$ in formula (\ref{eq:Oct27j1}) is quaternionic valued, while in this case it is complex valued (with values in $\mathbb C_J$).
One may wonder is there are choices of $I,J$ for which formula (\ref{an}) would allow to define $a_{-n}$ via (\ref{an}) and then obtain $a_{-n}=\bar a_n$. Since
$$
a_{-n}=2\int_0^{2\pi} e^{Int}d\mu_J(t),\qquad \qquad \bar a_n=2\int_0^{2\pi} \overline{d\mu_J(t)} e^{Int},
$$
and
\[
\begin{split}
\bar a_n&=2\int_0^{2\pi} (d\mu_1(t)-d\mu_2(t)J) e^{Int}\\
&=2\int_0^{2\pi} e^{Int} d\mu_1(t)-e^{-Int}d\mu_2(t)J \\
&=2\int_0^{2\pi} e^{Int} (d\mu_1(t)-d\mu_2(2\pi - t)J), \\
\end{split}
\]
 the condition $a_{-n}=\bar a_n$ translates into  ${\rm Re}(G)(e^ {It})=- {\rm Re}(G)(e^ {I(2\pi -t)})$. If one writes the power series expansion of $f_I$ in the form
 $$f_I(x+Iy)=\sum_{n=0}^\infty (x+Iy)^n a_n =\sum_{n=0}^\infty (x+Iy)^n(a_{0n}+Ia_{1n} +(a_{2n}+Ia_{3n})J)$$
 then it follows that
 $$
 F(x+Iy)=\sum_{n= 0}^\infty (x+Iy)^n(a_{0n}+Ia_{1n})
 \ \ \ {\rm and}\ \ \
   G(x+Iy)=\sum_{n = 0}^\infty (x+Iy)^n (a_{2n}+Ia_{3n}).
  $$
 Then $${\rm Re}(G)(x+Iy)=\sum_{n = 0}^\infty u_n(x,y) a_{2n} -v_n(x,y) a_{3n}$$
 where
 $$(x+Iy)^n=u_n(x,y)+Iv_n(x,y)$$
and
 $$
u_n(x,y)=\sum_{k=0,\, k\, even}^n {n\choose k}(-1)^{k/2}x^{n-k}y^k,
$$
$$
v_n(x,y)=\sum_{k=1,\, k\, odd}^n {n\choose k}(-1)^{(k-1)/2}x^{n-k}y^k.
$$
It is immediate that $u_n$ and $v_n$ are even and odd in the variable $y$, respectively, thus  ${\rm Re}(G)$ is odd in the variable $y$ if and only if $a_{2n}=0$ for all $n\in\mathbb N$. In general, given a slice hyperholomorphic function $f$ on $\mathbb B_1$ there is no change of basis for which one can have all the coefficients $a_{2n}=0$ for all $n\in\mathbb N$. Thus,  formula (\ref{an}) does not allow to define $a_{-n}$ in order to obtain the desired equality $a_{-n}=\bar a_n$. The formula is however one of the several possibilities to assign the coefficients of $f$ in integral form.
}
\end{Rk}
We conclude this section with a global integral representation.
\begin{Tm}
Let $f:\mathbb{B}_{1}\to \mathbb{H}$ be a slice hyperholomorphic function.
Let $f_I$ be the restriction of $f$ to the complex plane $\mathbb{C}_I$ and let
$F,G:  \mathbb{B}_{1}\cap \mathbb{C}_I \to \mathbb{C}_I$ be holomorphic functions $f_I(z)=F(z)+G(z)J$, $z=x+Iy$.
Then
$$
f(q)=I[{\rm Im} F(0) +{\rm Im} G(0)J] +\frac{1}{2\pi}\int_0^{2\pi}
K(q,e^{It}) \, d\mu_J(t),
$$
where $\mu_J(t)$ is a finite variation complex measure on $\mathbb{C}_J$  for $t\in [0,2\pi ] $ and
\[
\begin{split}
K(q,e^{It})&=
\frac 12\left(\frac{e^{It}+z}{e^{It}-z}+\frac{e^{It}+\overline{z}}{e^{It}-\overline{z}}\right)+\frac 12 I_qI
\left(\frac{e^{It}+\overline{z}}{e^{It}-\overline{z}} -\frac{e^{It}+z}{e^{It}-z}\right)
\\
&=
(1+q^2-2q{\rm Re}(e^{It}))^{-1}(1+2q{\rm Im}(e^{It})-q^2).
\end{split}
\]
\end{Tm}
\begin{proof}
From Theorem \ref{Herglotz} the restriction of $f$ to the complex plane $\mathbb{C}_I$ is
$$
f_I(z)=I [{\rm Im} F(0) +{\rm Im} G(0)J]+\int_0^{2\pi}\frac{e^{It}+z}{e^{It}-z} \, d\mu_J(t),
$$
where $\mu_J(t)$ is a finite variation complex measure on $\mathbb C_J$. Consider
$$
f_I(z)+f_I(\overline{z})
=2[{\rm Im} F(0) +{\rm Im} G(0)J]\left(\frac{e^{It}+z}{e^{It}-z}+\frac{e^{It}+\overline{z}}{e^{It}-\overline{z}}\right) \,d\mu_J(t)
$$
and
$$
f_I(\overline{z})-f_I(z)
=
\int_0^{2\pi}\left(\frac{e^{It}+\overline{z}}{e^{It}-\overline{z}} -\frac{e^{It}+z}{e^{It}-z}\right) \,d\mu_J(t);
$$
by applying the Representation Formula we obtain the kernel
$$
K(q,e^{It})=\frac 12\left(\frac{e^{It}+z}{e^{It}-z}+\frac{e^{It}+\overline{z}}{e^{It}-\overline{z}}\right)+\frac 12 I_qI
\left(\frac{e^{It}+\overline{z}}{e^{It}-\overline{z}} -\frac{e^{It}+z}{e^{It}-z}\right).
$$
 written in the first form. Now we write it in an equivalent way
 observing that the slice hyperholomorphic extension of the function
$$
K(z,e^{It})=\frac{e^{It}+z}{e^{It}-z}, \ \ \  z=x+Iy
$$
is (for the $\star$-inverse see Ch. 4 in \cite{MR2752913}) and the
$$
K(q,e^{It})=(e^{It}-q)^{-*}*(e^{It}+q)
$$
so that
$$
K(q,e^{It})=(1+q^2-2q{\rm Re}(e^{It}))^{-1}(1+2q{\rm Im}(e^{It})-q^2),
$$
and the statement follows.
\end{proof}

\bibliographystyle{plain}
\def\cprime{$'$} \def\lfhook#1{\setbox0=\hbox{#1}{\ooalign{\hidewidth
  \lower1.5ex\hbox{'}\hidewidth\crcr\unhbox0}}} \def\cprime{$'$}
  \def\cfgrv#1{\ifmmode\setbox7\hbox{$\accent"5E#1$}\else
  \setbox7\hbox{\accent"5E#1}\penalty 10000\relax\fi\raise 1\ht7
  \hbox{\lower1.05ex\hbox to 1\wd7{\hss\accent"12\hss}}\penalty 10000
  \hskip-1\wd7\penalty 10000\box7} \def\cprime{$'$} \def\cprime{$'$}
  \def\cprime{$'$} \def\cprime{$'$}

   \end{document}